\documentclass[12pt]{article}

\usepackage{amsmath}
\usepackage{theorem}
\usepackage{amssymb}
\usepackage{latexsym}
\usepackage{a4wide}

\newtheorem{thm}{Theorem}
\newtheorem{lem}[thm]{Lemma}

{\theorembodyfont{\rmfamily} }

\newenvironment{pf}{\noindent{\bf Proof.}}{\hbox{}\hfill $\Box$}

\newcommand{\Q}{\mathbb{Q}}
\newcommand{\R}{\mathbb{R}}
\newcommand{\C}{\mathbb{C}}
\newcommand{\F}{\mathbb{F}}

\newcommand{\Aut}{\mathrm{\mathop{Aut}}}
\newcommand{\GL}{\mathrm{\mathop{GL}}}

\begin{document}

\title{Classification of $6$-dimensional nilpotent Lie algebras over fields
of characteristic not $2$}
\author{Willem A. de Graaf\\
Dipartimento di Matematica\\
Universit\`{a} di Trento\\
Italy}
\date{}
\maketitle

\begin{abstract}
First we describe the Skjelbred-Sund method for classifying nilpotent 
Lie algebras. Then we use it to classify $6$-dimensional nilpotent Lie 
algebras over any field of characteristic not $2$. The proof of this 
classification is essentially constructive: for a given $6$-dimensional 
nilpotent Lie algebra $L$, following the steps of the proof, it is possible 
to find a Lie algebra $M$ that occurs in the list, 
and an isomorphism $L\to M$. 
\end{abstract}

\section{Introduction}

Several classifications of nilpotent Lie algebras of dimension
$6$ are available in the literature. We mention the ones by
V. Morozov (\cite{morozov}, over a field of characteristic $0$),
by R. Beck and B. Kolman (\cite{beko3}, over $\R$), by O. Nielsen
(\cite{nielsen}, over $\R$), by M.-P. Gong 
(\cite{gong}, over algebraically closed fields). Recently C. Schneider
has obtained lists of nilpotent Lie algebras over various finite
fields (\cite{schneider2}). \par
One of the main problems with these classifications is that they are
not easy to compare. One of the aims of this paper is to amend this
situation by describing a classification that is essentially algorithmic.
Starting from any nilpotent Lie algebra $L$, it is possible to follow the steps
in the proof of the classification to obtain an isomorphism between $L$ and
an element of the list. We call this the recognition procedure.
It is slightly tedious to perform this procedure by hand (see Section 
\ref{rec_exa}), but a computer can perform it very
efficiently. The recognition procedure has been implemented in the
computer algebra systems {\sf Magma} (\cite{Magma}) and {\sf GAP}4 
(\cite{GAP4}), the latter in the context of the {\sf liealgdb} package 
\cite{liealgdb}. It is expected that these implementations will be released 
in the near future. \par
The method that we use for classifying nilpotent Lie algebras of dimension $n$ 
essentially consists of two steps. In the first step a maybe redundant list of
Lie algebras is constructed that contains all $n$-dimensional nilpotent
Lie algebras. Then we erase isomorphic copies from the list. For the first
step we follow the method of T. Skjelbred and T. Sund (\cite{skjelsund}). The 
$n$-dimensional nilpotent Lie algebras are constructed as central extensions 
of nilpotent Lie algebras of smaller dimension. Subsequently the action
of the automorphism group is used to considerably reduce the number
of isomorphic Lie algebras occurring in the list. To get rid of the last
isomorphic copies, we use a technique for isomorphism testing based
on Gr\"{o}bner bases (cf. \cite{gra11}), which are computed in {\sf Magma}. 
This enables us to decide isomorphism, and
also to construct explicit isomorphisms. These are then 
used in the recognition procedure. We
remark that in many cases this method works very well. But in some
cases (e.g., if the Gr\"{o}bner basis is very complicated) we have not been 
able to use it. This is the reason why the paper does not contain the
classification over fields of characteristic $2$. \par
The results presented in this paper prove the following theorem.
\begin{thm}
Let $N_6(F)$ denote the number of $6$-dimensional nilpotent Lie algebras
over the field $F$. Suppose that the characteristic of $F$ is not $2$, and
let $s$ be the cardinality of $F^*/(F^*)^2$. Then $N_6(F)=26+4s$. In 
particular, $N_6(\C)=30$, $N_6(\R)=34$, $N_6(\Q)=\infty$,
$N_6(\F_q)=34$ for $q$ odd. 
\end{thm}

This paper is organised as follows. In Sections \ref{consect} to 
\ref{class_proc_sect} the Skjelbred-Sund method for classifying nilpotent
Lie algebras is described. We provide proofs of all results (absent
from \cite{skjelsund}) because some steps are needed in the recognition
procedure. Section \ref{rec_proc_sect} contains a sketch of this
procedure, and Section \ref{rec_exa} has an example. Then Section
\ref{sec_class_5} contains the classification of $5$-dimensional
nilpotent Lie algebras. And Section \ref{sec_class_6} has the
classification of $6$-dimensional nilpotent Lie algebras. Finally, in Section
\ref{sec_comm} we briefly comment on some of the classifications of
$6$-dimensional nilpotent Lie algebras that are found in the literature.\par
As mentioned, this paper does not contain a classification over 
fields of characteristic $2$. However, on some occasions we do comment
on what happens in that case.

\section{Constructing nilpotent Lie algebras}\label{consect}

Here and in the sequel we denote the centre of a Lie algebra $L$
by $C(L)$. The ground field of the vector spaces and Lie algebras 
will be denoted $F$. \par
Let $L$ be a Lie algebra, $V$ a vector space, and $\theta : L\times L
\to V$ a bilinear map. If $\theta(x,x)=0$ for all $x\in L$, then
$\theta$ is said to be skew-symmetric.\par
Let $\theta : L\times L\to V$ be a skew-symmetric bilinear map,
and set $L_{\theta} = L\oplus V$. For $x,y\in L$, $v,w\in V$ we define
$[x+v,y+w]=[x,y]_L+\theta(x,y)$. Then $L_{\theta}$ is a Lie algebra
if and only if 
$$\theta([x_1,x_2],x_3)+\theta([x_3,x_1],x_2)+\theta([x_2,x_3],x_1)=0
\text{ for all $x_1,x_2,x_3\in L$}.$$
The skew-symmetric $\theta$ satisfying this are called cocycles. The
set of all cocycles is denoted $Z^2(L,V)$. The Lie algebra $L_{\theta}$
is called a central extension of $L$ by $V$ (note that $V$ is central
in $L_{\theta}$). \par
Let $\nu : L\to V$ be a linear map, and define $\eta(x,y)=\nu([x,y])$.
Then $\eta$ is a cocycle, called a coboundary. The set of all coboundaries
is denoted $B^2(L,V)$. Let $\eta$ be a coboundary; then 
$L_{\theta}\cong L_{\theta+\eta}$.
Indeed, define $\sigma : L_{\theta}\to L_{\theta+\eta}$
by $\sigma(x) = x + \nu(x)$ for $x\in L$ and $\sigma(v)=v$ for $v\in V$.
The $\sigma$ is an isomorphism. Therefore we consider the set 
$H^2(L,V)=Z^2(L,V)/B^2(L,V)$. \par
Now let $K$ be a Lie algebra, and supppose that its centre, $C(K)$, is 
nonzero. Then we set $V=C(K)$, and $L=K/C(K)$. Let $\pi : K\to L$ be the
projection map. Choose an injective linear map $\sigma : L\to K$ such that
$\pi(\sigma(x))=x$ for all $x\in L$. Define $\theta : L\times L\to V$ by
$\theta(x,y) = [ \sigma(x), \sigma(y) ] -\sigma([x,y])$. 
Then $\theta$ is a cocycle. Note that $\theta$ depends on the choice of 
$\sigma$. But the $\theta$'s corresponding to two different $\sigma$'s differ
by a coboundary. Therefore, $\theta$ is a well-defined element of 
$H^2(L,V)$. Furthermore, $K\cong L_{\theta}$. Indeed, let $x\in K$,
then $x$ can uniquely be written as $x=\sigma(y)+z$, where $y\in L$ and
$z\in C(K)$. Define $\phi : K\to L_{\theta}$ by $\phi(x)=y+z$. Then $\phi$ 
is an isomorphism. (For that note that $[\sigma(y_1),\sigma(y_2)] =
\sigma([y_1,y_2])+\theta(y_1,y_2)$.)\par
We conclude that any Lie algebra with a nontrivial centre can be obtained
as a central extension of a Lie algebra of smaller dimension. So in
particular, all nilpotent Lie algebras can be constructed this way.

\section{The centre}

When constructing nilpotent Lie algebras as $L_{\theta}=L\oplus V$,
we want to restrict to $\theta$ such that $C(L_{\theta})=V$. If
the centre of $L_{\theta}$ is bigger, then it can be constructed as
a central extension of a different Lie algebra. (This way we avoid
constructing the same Lie algebra as central extension of different Lie
algebras.) There is a straightforward criterion on $\theta$ to decide
this. For $\theta\in Z^2(L,V)$ set
$$\theta^{\perp} =\{ x\in L\mid \theta(x,y)=0\text{ for all $y\in L$}\},$$
which is called the radical of $\theta$. Then $C(L_{\theta}) = 
(\theta^{\perp}\cap C(L))+V$, which immediately shows the following lemma.

\begin{lem}\label{lem2.1}
$\theta^{\perp}\cap C(L) =0$ if and only if $C(L_{\theta})=V$.
\end{lem}

\section{Isomorphism}

Let $e_1,\ldots,e_s$ be a basis of $V$, and $\theta \in Z^2(L,V)$. Then
$$\theta(x,y) = \sum_{i=1}^s \theta_i(x,y) e_i,$$
where $\theta_i \in Z^2(L,F)$. Furthermore, $\theta$ is
a coboundary if and only if all $\theta_i$ are.\par
Let $\phi\in \Aut(L)$. For $\eta\in Z^2(L,V)$ define $\phi\eta(x,y) = 
\eta(\phi(x),\phi(y))$. Then $\phi\eta \in Z^2(L,V)$. So $\Aut(L)$ acts on 
$Z^2(L,V)$. Also, $\eta\in B^2(L,V)$ if and only if $\phi\eta\in B^2(L,V)$
so $\Aut(L)$ acts on $H^2(L,V)$. 

\begin{lem}\label{lem1.4}
Let $\theta(x,y) = \sum_{i=1}^s \theta_i(x,y)e_i$ and 
$\eta(x,y)=\sum_{i=1}^s \eta_i(x,y) e_i$. Suppose that $\theta^{\perp}
\cap C(L)= \eta^{\perp}\cap C(L)=0$. Then
$L_{\theta}\cong L_{\eta}$ if and only if there is a $\phi\in\Aut(L)$
such that the $\phi\eta_i$ span the same subspace of $H^2(L,F)$
as the $\theta_i$.
\end{lem}

\begin{pf}
As vector spaces $L_{\theta} = L\oplus V$ and $L_{\eta} = L\oplus V$.
Let $\sigma : L_{\theta}\to L_{\eta}$ be an isomorphism. Since
$V$ is the centre of both Lie algebras, we have $\sigma(V) = V$. 
So $\sigma$ induces an isomorphism of $L_{\theta}/V=L$ to $L_{\eta}/V=L$,
i.e., an automorphim of $L$. Denote this automorphism by $\phi$. 
Let $L$ be spanned by $x_1,\ldots, x_n$. Then we write $\sigma(x_i) =
\phi(x_i)+v_i$, where $v_i\in V$, and $\sigma(e_i)=\sum_{j=1}^s a_{ji}e_j$.
Also write $[x_i,x_j]=\sum_{k=1}^n c_{ij}^k x_k$, and $v_i=\sum_{l=1}^s 
\beta_{il}e_l$. Then the relation $\sigma([x_i,x_j]_{L_{\theta}})=[\sigma(x_i),
\sigma(x_j)]_{L_{\eta}}$ amounts to
\begin{equation}\label{eq_isom} 
\eta_l(\phi(x_i),\phi(x_j)) = \sum_{k=1}^s a_{lk} \theta_k(x_i,x_j) 
+\sum_{k=1}^n c_{ij}^k \beta_{kl}, \text{ for $1\leq l\leq s$.}
\end{equation}
Now define the linear function $f_l : L\to F$ by $f_l(x_k) = \beta_{kl}$.
Then $f_l([x_i,x_j]) = \sum_k c_{ij}^k \beta_{kl}$. We see
that modulo $B^2(L,F)$, the $\phi \eta_i$ and the $\theta_i$ span the
same space.\par
Suppose that the $\phi \eta_i$ and $\theta_i$ span the same space
in $Z^2(L,F)$, modulo $B^2(L,F)$. Then there are linear functions 
$f_l : L\to F$ and $a_{lk}\in F$ so that $\phi\eta_l (x_i,x_j)
= \sum_k a_{lk}\theta_k (x_i,x_j) +f_l([x_i,x_j])$. If we set
$\beta_{kl}=f_l(x_k)$, then we see that (\ref{eq_isom}) holds.
This means that, if we define $\sigma : L_{\theta}\to L_{\eta}$ by
$\sigma(x_i)=\phi(x_i)+\sum_l \beta_{il}e_l$, $\sigma(e_i)=\sum_j a_{ji}
e_j$, then $\sigma$ is an isomorphism.
\end{pf}

\section{Avoiding central components}

Let $L=I_1\oplus I_2$ be the direct sum of two ideals. Suppose that 
$I_2$ is contained in the centre of $L$. Then $I_2$ is called a {\em central
component} of $L$. \par
The nilpotent Lie algebras with central components are simply obtained by 
taking direct sums of nilpotent Lie algebras of smaller dimension with 
abelian Lie algebras. Therefore, when constructing nilpotent Lie algebras
as central extensions we want to avoid constructing those with central 
components. For that we use the following criterion.

\begin{lem}\label{lem2.2}
Let $\theta$ be such that $\theta^{\perp}\cap C(L) =0$. Then
$L_{\theta}$ has no central components if and only if $\theta_1,
\dots,\theta_s$ are linearly independent in $H^2(L,F)$.
\end{lem}

\begin{pf}
Suppose that $L_{\theta}$ has a central component, $I$. Then we can write
$L_{\theta} = L\oplus W\oplus I$, where $L\oplus W \supset 
[L_{\theta},L_{\theta}]$.
Hence $\theta : L\times L\to W$, and we can write $\theta(x,y) =
\sum_{i=1}^t \tilde{\theta}_i (x,y) w_i$, where $t=\dim(W)$. By Lemma 
\ref{lem1.4}, the
$\tilde{\theta}_i$ span the same space as the $\theta_i$. It follows
that the $\theta_i$ are not linearly independent.\par
If the $\theta_i$ are linearly dependent, then by Lemma \ref{lem1.4} we may
assume that some of them are zero. This implies that $L_{\theta}$ has a 
central component. 
\end{pf}

\section{The classification procedure}\label{class_proc_sect}

Now we have a procedure that takes as input a nilpotent
Lie algebra $L$ of dimension $n-s$. It outputs all nilpotent Lie algebras $K$
of dimension $n$ such that $K/C(K)\cong L$, and $K$ has no central components.
It runs as follows.

\begin{enumerate}
\item Determine $Z^2(L,F)$ and $B^2(L,F)$ and determine the set 
$H^2(L,F)$ of cosets of $B^2(L,F)$ in $Z^2(L,F)$.
\item Consider $\theta \in H^2(L,V)$ with $\theta(x,y) = \sum_{i=1}^s 
\theta_i(x,y)e_i$ where
the $\theta_i\in H^2(L,F)$ are linearly independent, and $\theta^{\perp} 
\cap C(L)=0$. 
\item Find a (maybe redundant) list of representatives of the orbits of 
$\Aut(L)$ acting on the $\theta$ from 2.
\item For each $\theta$ found, construct $L_{\theta}$. Get rid of the 
isomorphic ones. 
\end{enumerate}

{\bf Comments:} For step 2. note that 
$$\theta^{\perp} = \theta_1^{\perp} \cap \cdots \cap \theta_s^{\perp}.$$
For step 3. there is no general method; this has to be done by hand,
on a case by case basis. \par
For step 4. we use Gr\"{o}bner bases to decide isomorphism, and
to construct isomorphisms (if they exist). This is described in detail 
in \cite{gra11}. There it is also mentioned that by writing some elements
of a Gr\"{o}bner basis of an ideal of a polynomial ring in terms of the 
polynomials that generate the ideal, we can decide isomorphism over
fields of all characteristics, apart from (maybe) a few exceptions. These
will then have to be considered separately. For nilpotent Lie algebras 
also a different approach is possible, based on Lemma \ref{lem1.4}. 
Suppose that we have the Lie algebras $L_{\theta}, L_{\eta}$, and we want
to decide whether they are isomorphic. Then we write a general element $\phi$ 
of $\Aut(L)$ by using indeterminates as its coordinates. Secondly, we 
introduce the indeterminates $\alpha_{ij}$ and we write the equations
$\phi\theta_i = \sum_j \alpha_{ij} \eta_j$. This leads to polynomial equations
in the coordinates of $\phi$ and the $\alpha_{ij}$. Finally, we decide 
solvability
by computing a Gr\"{o}bner basis of the ideal generated by these polynomials.
(This last part is completely analogous to the procedure described in 
\cite{gra11}.)
Examples of this approach are given in Sections \ref{sec_6_58}, 
\ref{sec_6_59}.

\section{The recognition procedure}\label{rec_proc_sect}

Let $K$ be a given nilpotent Lie algebra of dimension $\leq 6$. 
The proof of the classification given in this paper, yields an
algorithm for finding the Lie algebra $M$ occurring in the classification
that is isomorphic to $K$, as well as an isomorphism $K\to M$. \par
First we set $C=C(K)$ and $L=K/C$. Fix an injective linear map
$\sigma : L\to K$ such that $\pi\sigma(x)=x$ for all $x\in L$. Set
$\theta(x,y) = [\sigma(x),\sigma(y)]-\sigma([x,y])$. Then as seen in
Section \ref{consect}, we can explicitly construct an isomorphism
$\varphi : K\to L_{\theta}$.\par
Next we compute an isomorphism $\tau : L\to N$, where $N$ is a nilpotent
Lie algebra occurring in the classification. Let $\eta : N\times N\to C$ 
be defined by $\eta(x,y)=\theta(\tau^{-1}(x),\tau^{-1}(y))$. Then 
$\eta$ is a cocycle. Furthermore, $\tau$ extends to an isomorphism
$\tau : L_{\theta}\to N_{\eta}$ (by setting $\tau(u)=u$ for $u\in C$).\par
Write $\eta(x,y) = \sum_{i=1}^s \eta_i (x,y)e_i$, where the $e_i$ form
a basis of $C$. 
The proof of the classification provides an explicit method for
finding a $\phi\in \Aut(N)$ and $a_{ij}$ and a coclycle $\eta'$ 
and a coboundary $\eta''$ such that $N_{\eta'+\eta''}$ occurs in the 
classification, and $\phi\eta_i = \sum_j a_{ij}
\eta_j'$ (where $\eta'(x,y) = \sum_{i=1}^s \eta'_i(x,y) e_i$). Now the
construction in the proof of Lemma \ref{lem1.4} provides an isomorphism
$\nu : N_{\eta}\to N_{\eta'}$. Furthermore, in Section \ref{consect} it is 
shown how to constuct an isomorphism $\varepsilon : N_{\eta'}\to 
N_{\eta'+\eta''}$.\par
The final step consists of composing the isomorphisms:
$\varepsilon\circ\nu\circ\tau\circ\varphi : K\to N_{\eta'}$.

\section{Notation and useful facts}

Throughout the paper we use some notational conventions (which partly
follow \cite{gong}). If $L$ is a Lie algebra with basis 
$x_1,\ldots ,x_n$, then by $\Delta_{ij}$ we denote the skew-symmetric 
bilinear map $\Delta_{ij}: L\times L \to F$ with $\Delta_{ij}(x_i,x_j) =1$,
$\Delta_{ij}(x_j,x_i) =-1$, and $\Delta_{ij}(x_k,x_l) =0$ if 
$\{k,l\}\neq \{i,j\}$. On some occasions we also denote the characteristic
of the ground field by $\chi$. We use the following convention to describe 
the action of $\Aut(L)$ on $H^2(L,F)$. If the latter consists of e.g., 
$\theta_{a,b}=a\Delta_{13}+b\Delta_{24}$, and $\phi\in \Aut(L)$, then
$\phi\theta_{a,b} = \theta_{a',b'}$. Here $a',b'$ are expressed in terms
of $a,b$ and the coefficients of $\phi$. Then we write
$a\mapsto a'$ and $b\mapsto b'$. We describe the automorphism group of a
Lie algebra by giving the matrix of a general element. For this we use the
column convention: the action of a $\phi\in \Aut(K)$ on the $i$-th basis 
element of $K$ is given by the $i$-th column of the matrix of $\phi$.
The proof of the classification consists
of a number of subsections. The names of the Lie algebras (such as $K$, 
or $K^1$) are local to the subsection in which they occur, except when
the name is of the form $L_{d,k}$. The latter means that the correponding
Lie algebra occurs with the same name in the final list.

\begin{lem}\label{lem3.1}
Let $x,y\in F$ not both zero. Let $\delta\neq 0$ be given.
Then we can choose $a,b,c,d$ such that $ad-bc=\delta$ and
$$\begin{pmatrix} a & b \\
                  c & d \end{pmatrix} \begin{pmatrix} x\\
                                                      y \end{pmatrix}
= \begin{pmatrix} 1\\ 0 \end{pmatrix}.$$
\end{lem}

\begin{pf}
If $y=0$ then set $c=0$, $a=1/x$, $d=x\delta$. If $y\neq 0$ then set
$a=c=-y\delta$, $d=x\delta$, $b=(1+xy\delta)/y$.
\end{pf}

\begin{thm}\label{diagthm}
Let $\theta : L\times L\to F$ be a skew-symmetric bilinear form.
Then there is a basis of $L$ with respect to which $\theta =\Delta_{12}
+\Delta_{34}+\cdots +\Delta_{r,r+1}$. 
\end{thm}

A proof can for instance be found in \cite{jac2} (Theorem 7 of Chapter V).
We remark that the proof gives an explicit method for finding a basis
as in the theorem (or, equivalently, a linear transformation $\phi$
such that $\phi\theta$ has the described form).

\begin{lem}\label{lem3.2}
Let $a,b\in F$. Suppose that there are $x,y,s,t\in F$ such that 
$x^2a-y^2b=0$ and $s^2ab-t^2=0$ and $xy-st\neq 0$. Then there is an
$\alpha\in F^*$ such that $b=\alpha^2 a$.
\end{lem}

\begin{pf}
If $x$ and $y$ are both nonzero, then this follows from the first
equation. If $x=0$, then $y=0$ or $b=0$. However, $b=0$ implies $t=0$
and $xy-st=0$. So $y=0$, and hence $s,t$ are both nonzero. But then
$a= (\frac{t}{bs})^2b$. If $y=0$ then we proceed in a similar way. 
\end{pf}

There are the following nilpotent Lie algebras of dimensions $3,4$
(cf. e.g., \cite{gra11}):

\begin{enumerate}
\item[$L_{3,1}$] Abelian.
\item[$L_{3,2}$] $[x_1,x_2]=x_3$.
\item[$L_{4,1}$] Abelian.
\item[$L_{4,2}$] $[x_1,x_2]=x_3$.
\item[$L_{4,3}$] $[x_1,x_2]=x_3, [x_1,x_3]=x_4.$
\end{enumerate}

\section{Nilpotent Lie algebras of dimension 5}\label{sec_class_5}

Here we use the method to classify nilpotent Lie algebras of dimension $5$.
We have a number of subsections; in each subsection the central extensions
of one particular Lie algebra are considered. 

\subsection{$L_{4,1}$}

In this case $H^2(L_{4,1},F)$ consists of 
all skew-symetric $\theta : L_{4,1}\times L_{4,1}\to F$. 
Let $\theta\in H^2(L_{4,1},F)$ be such that $\theta^{\perp}
\cap C(L_{4,1}) = 0$. This means that the matrix of $\theta$ is non-singular.
By Theorem \ref{diagthm},
there is a basis of $L_{4,1}$ with respect to which $\theta = \Delta_{1,2}+
\Delta_{3,4}$. We get the Lie algebra
$$L_{5,4}: [x_1,x_2] = x_5, [x_3,x_4]=x_5.$$

\subsection{$L_{4,2}$}\label{secnonab4_1}

Here $H^2(L_{4,2},F)$ consists of all $a\Delta_{13}+b\Delta_{14}+
c\Delta_{23}+d\Delta_{24}$. $\Aut(L_{4,2})$ consists of 
$$\phi= \begin{pmatrix}  a_{11} & a_{12} & 0 & 0 \\
                         a_{21} & a_{22} & 0 & 0 \\
                         a_{31} & a_{32} & \delta & a_{34}\\
                         a_{41} & a_{42} & 0 & a_{44} \end{pmatrix},$$
where $\delta = a_{11}a_{22}-a_{12}a_{21}\neq 0$. 
The automorphism $\phi$ as above acts as follows
\begin{align*}
a &\mapsto \delta( a_{11}a+a_{21}c)\\
b &\mapsto a_{11}a_{34}a + a_{11}a_{44}b + a_{21}a_{34}c+a_{21}a_{44}d\\
c &\mapsto \delta( a_{12}a+a_{22}c )\\
d &\mapsto a_{12}a_{34}a+a_{12}a_{44}b+a_{22}a_{34}c+a_{22}a_{44}d.
\end{align*}

Here we determine all orbits of $H^2(L_{4,2},F)$ under the action of 
$\Aut(L_{4,2})$, because we need that later. \par
First suppose that $a\neq 0$. Then we set $a_{11}=1/a$, $a_{22}=a$, and
$a_{21}=0$. Then $a\mapsto 1$. To conserve this situation we set 
$a_{21}=0$ and $a_{11}=a_{22}=1$. With $a_{12}=-c$ we get 
$c\mapsto 0$. To conserve this we set $a_{12}=0$. Then $b\mapsto 
a_{34}+a_{44}b$. So we set $a_{34}=-a_{44}b$, and $b\mapsto 0$.
Also, $d\mapsto a_{44} d$. Therefore, if $d\neq 0$ we can get
$d\mapsto 1$. Hence we get two cocycles (corresponding to $d=0,1$):
$\theta=\Delta_{13}$ and $\theta=\Delta_{13}+\Delta_{24}$. \par
If $a=0$, but $c\neq 0$, then we choose $a_{21}\neq 0$ and
get that $a$ is mapped to something nonzero. Hence we are in the previous
case.\par
Suppose that $a=c=0$. Then $b\mapsto a_{44}(a_{11}b+a_{21}d)$, and
$d\mapsto a_{44}(a_{12}b+a_{22}d)$. So by Lemma \ref{lem3.1}, 
we see that we can choose tha $a_{ij}$ in such a way that $b$ is
mapped to $1$, and $d$ is mapped to $0$. Hence we get $\theta=\Delta_{14}$.\par
We have obtained three cocycles. Any two of them are not in the same orbit, as
$\Delta_{13}$ is the only one in the subspace spanned by 
$\Delta_{13}$ and $\Delta_{23}$, which is stable under $\Aut(L)$.
And the other two have a radical of a different dimension.\par
The centre of $L_{4,2}$ is spanned by $x_3,x_4$. So only one of the
cocycles obtained has a radical which has zero intersection
with $C(L_{4,2})$, namely   $\Delta_{13}+\Delta_{24}$. Hence we get 
one Lie algebra:
$$L_{5,5}: [x_1,x_2]=x_3, [x_1,x_3]= x_5, [x_2,x_4] = x_5.$$

\subsection{$L_{4,3}$}

$H^2(L_{4,3},F)$ consists of $a\Delta_{14}+b\Delta_{23}$. The centre of 
$L_{4,3}$ is spanned by $x_4$. So $\theta^{\perp}\cap C(L_{4,3})=0$ if and 
only if $a\neq 0$. Hence we may normalise to $a=1$. We get the Lie algebras
$$[x_1,x_2]=x_3, [x_1,x_3]=x_4, [x_1,x_4]=x_5, [x_2,x_3]=bx_5.$$
If $b\neq 0$ we multiply $x_1,\ldots, x_5$ by $b,b,b^2,b^3,b^4$, and we
see that this Lie algebra is isomorphic to the one with $b=1$. Hence
we get two Lie algebras: $L_{5,6}$ ($b=1$) and $L_{5,7}$ ($b=0$). 
A Gr\"{o}bner basis computation shows that they are not isomorphic.

\subsection{$L_{3,1}$}

Here $H^2(L_{3,1},F)$ consists of $a\Delta_{12}+b\Delta_{13}+c\Delta_{23}$.
If $\theta$ is of this form, then $\theta^{\perp}$ is spanned by 
$cx_1-bx_2+ax_3$. Hence if we take two linearly independent elements
of $H^2(L_{3,1},F)$ then the intersection of the radicals is zero.\par
An element $\phi = (a_{ij})\in \Aut(L_{3,1})$ acts as follows on
$H^2(L_{3,1},F)$:
\begin{align*}
a &\mapsto (a_{11}a_{22}-a_{12}a_{21})a+ (a_{11}a_{32}-a_{31}a_{12})b+
(a_{21}a_{32}-a_{31}a_{22})c\\
b &\mapsto (a_{11}a_{23}-a_{21}a_{13})a + (a_{11}a_{33}-a_{31}a_{13})b +
(a_{21}a_{33}-a_{31}a_{23})c \\
c &\mapsto (a_{12}a_{23}-a_{22}a_{13})a+(a_{12}a_{33}-a_{32}a_{13})b +
(a_{22}a_{33}-a_{32}a_{23})c.
\end{align*}

Here we classify $2$-dimensional subspaces of $H^2(L_{3,1},F)$. Let such a 
subspace
by spanned by $\theta_1$, $\theta_2$. By Theorem \ref{diagthm} we may assume
that $\theta_1 = \Delta_{12}$. In the sequel we consider automorphisms
$\phi$ with $a_{23}=a_{13}=0$. This means that $\phi\theta_1$ is a scalar
multiple of $\theta_1$. After subtracting a suitable multiple of $\theta_1$
we may assume that $\theta_2 = b\Delta_{13}+c\Delta_{13}$. We can choose 
$\phi$ such that $b\mapsto b'\neq 0$, hence we may assume that $b\neq 0$.
Then we divide and get $b=1$. Now under $\phi$, $b\mapsto a_{11}a_{33}b
+a_{21}a_{33}c$, $c\mapsto a_{12}a_{33}b+a_{22}a_{33}c$. So if we set
$a_{11}=a_{22}=a_{33}=1$, $a_{21}=0$ and $a_{12}=-c$, then $c\mapsto 0$. 
Hence we get $\theta_2=\Delta_{13}$ and the Lie algebra:
$$L_{5,8}: [x_1,x_2]=x_4, [x_1,x_3]=x_5.$$

\subsection{$L_{3,2}$}

$H^2(L_{3,2},F)$ is 
two dimensional, and spanned by $\Delta_{13},\Delta_{23}$. So there is
only one subspace of dimension $2$. This leads to the Lie algebra
$$L_{5,9}: [x_1,x_2]=x_3, [x_1,x_3]=x_4, [x_2,x_3]=x_5.$$

\subsection{The list}

We get nine nilpotent Lie algebras of dimension $5$. First we have
$L_{5,k}$ for $k=1,2,3$ that are the direct sum of $L_{4,k}$ with a 
$1$-dimensional abelian ideal. The $5$-dimensional nilpotent Lie algebras
without central component are:
\begin{enumerate}
\item[$L_{5,4}$] $[x_1,x_2] = x_5, [x_3,x_4]=x_5.$
\item[$L_{5,5}$] $[x_1,x_2]=x_3, [x_1,x_3]= x_5, [x_2,x_4] = x_5.$
\item[$L_{5,6}$] $[x_1,x_2]=x_3, [x_1,x_3]=x_4, [x_1,x_4]=x_5, [x_2,x_3]=x_5.$
\item[$L_{5,7}$] $[x_1,x_2]=x_3, [x_1,x_3]=x_4, [x_1,x_4]=x_5.$
\item[$L_{5,8}$] $[x_1,x_2]=x_4, [x_1,x_3]=x_5.$
\item[$L_{5,9}$] $ [x_1,x_2]=x_3, [x_1,x_3]=x_4, [x_2,x_3]=x_5.$
\end{enumerate}

\section{Nilpotent Lie algebras of dimension 6}\label{sec_class_6}

We list central extensions of the nilpotent Lie algebras of 
dimension $\leq 5$.\par
Here we do not need to consider central extensions of $L_{5,1}$, as
there are no non-degenerate skew-symmetric bilinear forms on 
an odd-dimensional space (Theorem \ref{diagthm}).

\subsection{$L_{5,2}$}

The automorphism group of $L_{5,2}$ consists of
 
$$ \begin{pmatrix} 
a_{11} & a_{12} & 0 & 0 & 0 \\
a_{21} & a_{22} & 0 & 0 & 0\\
a_{31} & a_{32} & a_{11}a_{22}-a_{12}a_{21} & a_{34} & a_{35} \\
a_{41} & a_{42} & 0 & a_{44} & a_{45}\\
a_{51} & a_{52} & 0 & a_{54} & a_{55}
\end{pmatrix}.$$
The elements of $H^2(L_{5,2},F)$ are $\theta=a\Delta_{13}+b(\Delta_{14}+
\Delta_{15})+
c\Delta_{15}+d\Delta_{23}+e\Delta_{24}+f\Delta_{25}+g\Delta_{45}$, and 
$\theta^{\perp}\cap C(L_{5,2})=0$ if and only if $g\neq 0$ and one of $a,d$ is 
not $0$. The group acts as follows (we put $\delta=a_{11}a_{22}-a_{12}a_{21}$)
\begin{align*}
a & \mapsto \delta(a_{11}a+a_{21}d)\\
b & \mapsto a_{11}a_{34}a+a_{11}a_{44}b+a_{11}a_{54}c+a_{21}a_{34}d+a_{21}
a_{44}e+a_{21}a_{54}f +(a_{41}a_{54}-a_{51}a_{44})g \\
c & \mapsto a_{11}a_{35}a+a_{11}a_{45}b+a_{11}a_{55}c+a_{21}a_{35}d+a_{21}
a_{45}e+a_{21}a_{55}f +(a_{41}a_{55}-a_{51}a_{45})g \\
d & \mapsto \delta(a_{12}a+a_{22}d)\\
e & \mapsto a_{12}a_{34}a+a_{12}a_{44}b+a_{12}a_{54}c+a_{22}a_{34}d+a_{22}
a_{44}e+a_{22}a_{54}f +(a_{42}a_{54}-a_{52}a_{44})g \\
f & \mapsto a_{12}a_{35}a+a_{12}a_{45}b+a_{12}a_{55}c+a_{22}a_{35}d+a_{22}
a_{45}e+a_{22}a_{55}f +(a_{42}a_{55}-a_{52}a_{45})g \\
g & \mapsto (a_{44}a_{55}-a_{45}a_{54})g.
\end{align*}

We map $\theta$ to a ``normal form'' in stages. Since not both $a$ and
$d$ are zero, we can choose $a_{11}, a_{12}, a_{21}, a_{22}$ such that
$a_{11}a+a_{21}d=1$, $a_{12}a+a_{22}d=0$ and $a_{11}a_{22}-a_{12}a_{21}=1$
(Lemma \ref{lem3.1}). Then $a\mapsto 1$ and $d\mapsto 0$.\par
Now we assume that we have $\theta$ with $a=1$ and $d=0$. Then we choose
$a_{11}=a_{22}=a_{44}=a_{55}=1$ and $a_{21}=a_{12}=a_{41}=a_{51}=
a_{54}=a_{45}=0$ and $a_{34}=-b$, $a_{35}=-c$. 
Then $a,b,c,d$ are mapped to $1,0,0,0$. \par
Now we assume that $a=1$, $b=c=d=0$. To conserve this we set 
$a_{11}=a_{22}=1$ and $a_{12}=a_{21}=a_{34}=a_{35}=a_{41}=a_{51}=a_{42}=0$. 
Now with $a_{44}=a_{55}=1$ and $a_{54}=a_{45}=0$ we get 
Then $e\mapsto e - a_{52}g$, $f\mapsto f+a_{42}g$. So since $g\neq 0$, we
can choose $a_{42},a_{52}$ such that $e,f$ are maped to $0$. \par 
Hence we may assume that $a=1$ and $b=c=d=e=f=0$. 
We put $a_{44}=1/g$, $a_{ii}=1$ for $i\neq 4$ and $a_{ij}=0$ otherwise.
Then $g\mapsto 1$ and we get the Lie algebra
$$L_{6,10} : [x_1,x_2]=x_3, [x_1,x_3]=x_6, [x_4,x_5]=x_6.$$

\subsection{$L_{5,3}$}

The automorphism group of $L_{5,3}$ consists of
 
$$ \begin{pmatrix} a_{11} & 0 & 0 & 0 & 0 \\
                   a_{21} & a_{22} & 0 & 0 & 0\\
                   a_{31} & a_{32} & a_{11}a_{22} & 0 & 0 \\
                   a_{41} & a_{42} & a_{11}a_{32} & a_{11}^2a_{22} & a_{45}\\
                   a_{51} & a_{52} & 0 & 0 & a_{55}
\end{pmatrix}.$$
The elements of $H^2(L_{5,3},F)$ are $\theta=a\Delta_{14}+b\Delta_{15}+
c\Delta_{23}+d\Delta_{25}$, and $\theta^{\perp}\cap C(L_{5,3})=0$ if and 
only if both $a$ and $d$ are nonzero. The group acts as follows
\begin{align*}
a & \mapsto a_{11}^3a_{22}a\\
b & \mapsto a_{11}a_{45}a+a_{11}a_{55}b+a_{21}a_{55}d\\
c & \mapsto a_{11}a_{22}^2c\\
d & \mapsto a_{22}a_{55}d. 
\end{align*}

After dividing we get $a=1$. Then we take $a_{11}=a_{22}=1$, $a_{55}=
1/d$, $a_{45}=0$ and $a_{21}=-b/d$. This means that 
$a,b,c,d$ are mapped to $1,0,c,1$ respectively. The corresponding Lie
algebra is given by 
$$K_c: [x_1,x_2]=x_3, [x_1,x_3]=x_4, [x_1,x_4]=x_6, [x_2,x_3]=cx_6,
[x_2,x_5]=x_6.$$
If $c\neq 0$ then we multiply $x_1,\ldots, x_6$ respectively by 
$c,c,c^2,c^3,c^3,c^4$, and get that $K_c\cong K_1$. A Gr\"{o}bner basis
computation shows that $K_1$ and $K_0$ are not isomorphic. So we get two
algebras: $L_{6,11}$ ($=K_1$) and $L_{6,12}$ ($=K_0$).

\subsection{$L_{5,4}$}

In this case $H^2(L_{5,4},F)$ contains the elements $\theta = a\Delta_{13}
+b\Delta_{14}+c\Delta_{23}+d\Delta_{24}+e\Delta_{34}$. However, each
such $\theta$ has $x_5$ in its radical. Hence there are no central
extensions of $L_{5,4}$ with a $1$-dimensional centre.

\subsection{$L_{5,5}$}

The automorphism group of $L_{5,5}$ consists of 
 
$$ \begin{pmatrix} a_{11} & 0 & 0 & 0 & 0 \\
                   a_{21} & a_{22} & 0 & 0 & 0\\
                   a_{31} & a_{32} & a_{11}a_{22} & -a_{11}a_{21} & 0 \\
                   a_{41} & a_{42} & 0 & a_{11}^2 & 0 \\
                   a_{51} & a_{52} & \varepsilon & a_{54} & a_{11}^2a_{22}
\end{pmatrix},$$
with $\varepsilon = a_{11}a_{32}+a_{21}a_{42}-a_{41}a_{22}$. \par
Here $H^2( L_{5,5}, F)$ consists of the elements $\theta = a\Delta_{14}
+b(\Delta_{15}+\Delta_{34})+c\Delta_{23}+d\Delta_{24}$. In order to ensure
that $\theta$ does not have $x_5$ (which spans $C(L_{5,5})$) in its radical, 
we must have $b\neq 0$. The automorphism group acts as follows:
\begin{align*}
a & \mapsto a_{11}^3a + (a_{11}a_{54}+a_{31}a_{11}^2+a_{11}a_{41}a_{21})b
-a_{11}a_{21}^2c +a_{21}a_{11}^2 d\\
b & \mapsto a_{11}^3a_{22}b\\
c & \mapsto a_{11}a_{22}^2 c -a_{11}a_{22}a_{42}b\\
d & \mapsto (a_{32}a_{11}^2+a_{42}a_{11}a_{21})b -a_{11}a_{22}a_{21}c+
a_{11}^2a_{22}d.
\end{align*}

Since $b\neq 0$ after dividing we may assume that $b=1$. Then we set
$a_{11}=a_{22}=1$, $a_{54}=a_{21}=0$, $a_{31}=-a$, $a_{42}=c$. This leads
to $a=c=0$. This yields the Lie algebras
$$K_d: [x_1,x_2]=x_3, [x_1,x_3]=x_5, [x_2,x_4]=x_5+dx_6, [x_1,x_5]=x_6,
[x_3,x_4]=x_6.$$
If $d\neq 0$ then we multiply $x_1,\ldots, x_6$ by $d,1,d,d^2,d^2,d^3$,
and we see that $K_d\cong K_1$. If the characteristic is not $2$, then 
$K_1\cong K_0$ by $x_1\mapsto y_1+\frac{1}{2}y_4$, $x_2\mapsto 
y_2+y_3+\frac{1}{2}y_5$, $x_3\mapsto y_3+\frac{1}{2}y_5$, $x_4\mapsto y_4$,
$x_5\mapsto y_5$, $x_6\mapsto y_6$. Hence we get one algebra, $L_{6,13}=K_0$.
\par
If the characteristic is $2$, then
a Gr\"{o}bner basis computation shows that $K_0$ and $K_1$ are not
isomorphic. So in that case we get two algebras.

\subsection{$L_{5,6}$}

The automorphism group consists of
$$ \begin{pmatrix} a_{11} & 0 & 0 & 0 & 0 \\
                   a_{21} & a_{11}^2 & 0 & 0 & 0\\
                   a_{31} & a_{32} & a_{11}^3 & 0 & 0 \\
                   a_{41} & a_{42} & a_{11}a_{32} & a_{11}^4 & 0 \\
                   a_{51} & a_{52} & u & v & a_{11}^5
\end{pmatrix},$$
where $u=a_{11}a_{42}+a_{21}a_{32}-a_{31}a_{11}^2$, $v=a_{21}a_{11}^3
+a_{32}a_{11}^2$. \par
As representatives of the elements of $H^2(L_{5,6},F)$ we use 
$\theta = a(\Delta_{15}+\Delta_{24})+b\Delta_{23}+c(\Delta_{25}-\Delta_{34})$. 
Moreover, $\theta^{\perp}\cap C(L_{5,6})=0$ if and only if $a\neq 0$ or 
$c\neq 0$. 
In order to compute the action of the automorphism group, some care is needed.
Let $\theta\in Z^2(L_{5,6},F)$ be given by $\theta = \alpha\Delta_{1,2}
+\beta\Delta_{13}+\gamma(\Delta_{14}+\Delta_{23})+a(\Delta_{15}+\Delta_{24})
+b\Delta_{23}+c(\Delta_{25}-\Delta_{34})$. Let the coefficients
of $\phi\theta$ be denoted by $\alpha', \beta', \ldots, c'$. Then
\begin{align*}
a' &= a_{11}^5(a_{11}a+a_{21}c)\\
b' &= -2a_{11}^4a_{21}a+a_{11}^5b+(2a_{11}^3a_{42}-a_{11}^3a_{21}^2
-a_{11}a_{32}^2)c\\
c' &= a_{11}^7 c.
\end{align*}
We distinguish two cases. If $c\neq 0$, then we divide and get $c=1$.
Then if the characteristic is not $2$ we choose $a_{11}=1$, $a_{32}=0$, 
$a_{21}=-a$, $a_{42}=
-\frac{1}{2}(a^2+b)$. This yields $a=b=0$. We get the Lie algebra
$$L_{6,14}: [x_1,x_2]=x_3, [x_1,x_3]=x_4, [x_1,x_4]=x_5, [x_2,x_3]=x_5, 
[x_2,x_5]=x_6, [x_3,x_4]=-x_6.$$
If $\chi=2$, then we choose $a_{32}$ so that $a_{32}^2=a^2+b$. If the field
is perfect this can be done. This leads to the same result.\par
If $c=0$ then $a\neq 0$ and we divide and get $a=1$. If $\chi\neq 2$ then
with $a_{11}=1$ and $a_{21}=b/2$ we get $b=0$. The corresponding Lie algebra is
$$L_{6,15}: [x_1,x_2]=x_3, [x_1,x_3]=x_4, [x_1,x_4]=x_5, [x_2,x_3]=x_5, 
[x_1,x_5]=x_6, [x_2,x_4]=x_6.$$
A Gr\"{o}bner basis computation shows that these two Lie algebras are not
isomorphic. If $\chi=2$ then we get a parametrized series of Lie algebras.

\subsection{$L_{5,7}$}

The automorphism group consists of
$$ \begin{pmatrix} a_{11} & 0 & 0 & 0 & 0 \\
                   a_{21} & a_{22} & 0 & 0 & 0\\
                   a_{31} & a_{32} & a_{11}a_{22} & 0 & 0 \\
                   a_{41} & a_{42} & a_{11}a_{32} & a_{11}^2a_{22} & 0 \\
                   a_{51} & a_{52} & a_{11}a_{42} & a_{11}^2a_{32} & 
a_{11}^3a_{22}
\end{pmatrix}.$$
$H^2(L_{5,7},F)$ consists of $\theta = a\Delta_{15}+b\Delta_{23}+c(\Delta_{25}
-\Delta_{34})$. Moreover, $\theta^{\perp}\cap C(L_{5,7})=0$ if and
only if $a\neq 0$ or $c\neq 0$. The automorphism group acts in the following 
way:
\begin{align*}
a & \mapsto a_{11}^3(a_{11}a_{22}a+a_{22}a_{21}c)\\
b & \mapsto a_{11}(a_{22}^2b + (2a_{22}a_{42}-a_{32}^2)c)\\
c & \mapsto a_{11}^3a_{22}^2 c.
\end{align*}

If $c\neq 0$ then we divide to get $c=1$. Choose $a_{11}=a_{22}=1$ and
$a_{21}=-a$, which yields $a=0$. If the characteristic of the ground field
is not $2$, then we set $a_{42}=-\frac{1}{2}b$ and $a_{32}=0$, and we
get $b=0$. (We remark
that, if $F$ is a perfect field of characteristic $2$, then we can
choose $a_{32}=\beta$ with $\beta^2=b$, leading to the same result.)
The corresponding Lie algebra is
$$L_{6,16}: [x_1,x_2]=x_3, [x_1,x_3]=x_4, [x_1,x_4]=x_5, [x_2,x_5]=x_6,
[x_3,x_4]=-x_6.$$
If $c=0$ then $a\neq 0$, and after dividing, $a=1$. We get the Lie algebras
$$K_b: [x_1,x_2]=x_3, [x_1,x_3]=x_4, [x_1,x_4]=x_5, [x_1,x_5]=x_6,
[x_2,x_3]=bx_6.$$
If $b\neq 0$ then we multiply $x_i$ by $b^i$ and we see that $K_b\cong K_1$.
A Gr\"{o}bner basis computation shows that $K_0\not\cong K_1$, and these
two algebras are also not isomorphic to $L_{6,16}$. So we get
two more algebras: $L_{6,17}$ ($=K_1$), and $L_{6,18}$ ($=K_0$).

\subsection{$L_{5,8}$}\label{sec_6_58}

The automorphism group consists of
\begin{equation}\label{autL58}
 \begin{pmatrix} a_{11} & 0 & 0 & 0 & 0 \\
                   a_{21} & a_{22} & a_{23} & 0 & 0\\
                   a_{31} & a_{32} & a_{33} & 0 & 0 \\
                   a_{41} & a_{42} & a_{43} & a_{11}a_{22} & a_{11}a_{23} \\
                   a_{51} & a_{52} & a_{53} & a_{11}a_{32} & a_{11}a_{33}
\end{pmatrix}.
\end{equation}

$H^2(L_{5,8},F)$ consists of $\theta = a\Delta_{14}+b\Delta_{15}+c\Delta_{23}
+d\Delta_{24}+e(\Delta_{25}+\Delta_{34})+f\Delta_{35}$. Moreover, 
$\theta^{\perp}\cap C(L_{5,8})=0$ if and
only if the matrix 
\begin{equation}\label{eqnL58}
\begin{pmatrix} a & b \\ d & e \\ e & f \end{pmatrix}
\end{equation}
has rank $2$. The automorphism group acts in the following 
way:
\begin{align*}
a & \mapsto a_{11}^2(a_{22}a+a_{32}b)+a_{11}a_{22}a_{21}d+a_{11}(a_{32}
a_{21}+a_{31}a_{22})e+a_{11}a_{32}a_{31}f\\
b & \mapsto a_{11}^2(a_{23}a+a_{33}b)+a_{21}a_{11}a_{23}d+a_{11}(a_{21}a_{33}
+a_{31}a_{23})e +a_{11}a_{32}a_{31}f\\
c & \mapsto a_{22}a_{33}c+(a_{22}a_{43}-a_{42}a_{23})d+
(a_{22}a_{53}+a_{32}a_{43}-a_{42}a_{33}-a_{52}a_{23})e+(a_{32}a_{53}-a_{52}
a_{33})f\\
d & \mapsto a_{11}a_{22}^2 d +2a_{11}a_{22}a_{32}e + a_{32}^2a_{11}f\\
e & \mapsto a_{11}a_{22}a_{23}d + a_{11}(a_{22}a_{33}+a_{23}a_{32})e
+a_{11}a_{33}a_{32}f\\
f & \mapsto a_{11}a_{23}^2 d + 2a_{11}a_{33}a_{23} e + a_{11}a_{33}^2 f. 
\end{align*}

First suppose $e\neq 0$. 
Note that $e\mapsto a_{23}(a_{11}a_{22}d+a_{11}a_{32}e)+a_{11}a_{33}
(a_{22}e+a_{32}f)$. We claim that we can map $e$ to $0$ unless
$d=f=0$ and the characteristic of the ground field is $2$. If $f\neq 0$,
then choose $a_{23}=0$ and $a_{32}=-a_{22}e/f$. If $f=0$, but $d\neq 0$
then choose $a_{32}=0$ and $a_{23}=-a_{33}e/d$. If $d=0$ as well, then 
$e\mapsto a_{11}(a_{33}a_{22}+a_{23}a_{32})e$. If the 
characteristic of the ground field is not $2$, then choose 
$a_{23}$ and $a_{32}$ in such a way that $a_{22}a_{33}-a_{23}a_{32}\neq 0$
but $a_{23}a_{32}+a_{22}a_{33}=0$. So if not $d=f=0$ and $\chi=2$,
then we may assume that 
$e=0$. To preserve this situation we set $a_{23}=a_{32}=0$. We distinguish
two cases. 
\begin{enumerate}

\item $d\neq 0$. Then we divide to get $d=1$. To preserve this we
choose $a_{11}=a_{22}=1$. Then choose $a_{21}$ and
$a_{43}$ so that $a\mapsto 0$ and $c\mapsto 0$. Again we consider
two cases.
\begin{enumerate}
\item $b=0$. Then we get the Lie algebras
$$L_{6,19}(\epsilon):  [x_1,x_2]=x_4, [x_1,x_3]=x_5, [x_2,x_4]=x_6, 
[x_3,x_5]=\epsilon x_6.$$
(Here we write $\epsilon$ in place of $f$. We have $\epsilon\neq 0$, 
as otherwise the matrix (\ref{eqnL58}) does not have rank $2$,
or in other words, the centre has dimension $>1$.)
Setting $y_i=x_i$ for $i=1,2,4,6$ and $y_i=\alpha x_i$ otherwise,
we see that $L_{6,19}(\epsilon)\cong L_{6,19}(\alpha^2 \epsilon)$. 
Now suppose that $L_{6,19}(\epsilon)\cong L_{6,19}(\delta)$. These algebras 
are defined relative to the cocycles
$\theta_{\epsilon}=\Delta_{24}+\epsilon\Delta_{35}$ and 
$\theta_{\delta}=\Delta_{24}+\delta
\Delta_{35}$. So, according to Lemma \ref{lem1.4}, there is a 
$\phi\in \Aut(L_{5,8})$ such that $\phi \theta_{\epsilon} = \alpha 
\theta_{\delta}$. Now
with $\phi$ as in (\ref{autL58}), we have that $\phi \theta_{\epsilon}=  
a_1\Delta_{14}+b_1\Delta_{15}+c_1\Delta_{23}
+d_1\Delta_{24}+e_1(\Delta_{25}+\Delta_{34})+f_1\Delta_{35}$, where
$a_1=a_{11}(a_{22}a_{21}+a_{32}a_{31}\epsilon)$, $b_1 = a_{11}(a_{11}a_{23}+
a_{32}a_{31}\epsilon)$, $c_1=a_{22}a_{43}-a_{42}a_{23}+(a_{32}a_{53}-
a_{52}a_{33})\epsilon$,
$d_1=a_{11}(a_{22}^2+a_{32}^2\epsilon)$, $e_1=a_{11}(a_{22}a_{23}+a_{33}
a_{32}\epsilon)$, $f_1=a_{11}(a_{23}^2+a_{33}^2\epsilon)$. Since 
$a_1,b_1,c_1,e_1$ have to be $0$
we get some polynomial equations (where we factor out the $a_{11}$ since 
it cannot be $0$). Also, $d_1=\alpha$ and $f_1=\alpha \delta$. This leads to
the equation $a_{23}^2+a_{33}^2\epsilon -(a_{22}^2+a_{32}^2\epsilon)\delta=0$. 
Finally we add
two equations that express that the matrix of the automorphism is nonzero:
$D_1a_{11}-1=0$, $D_2( a_{22}a_{33}-a_{32}a_{23})-1=0$ and to force
$d_1=\alpha\neq 0$: $D_3(a_{22}^2+a_{32}^2\epsilon)-1=0$.
The Gr\"{o}bner basis of this system contains $a_{22}^2\delta - 
a_{33}^2\epsilon$ and $a_{32}^2\epsilon\delta - a_{23}^2$ (over any field). 
So by Lemma \ref{lem3.2} 
we conclude that $\epsilon=\beta^2\delta$ for some $\beta\in F^*$. 
\item $b\neq 0$. Then choose $a_{33}= 1/b$ to get 
$b=1$. We get the Lie algebras
$$K^1_f:  [x_1,x_2]=x_4, [x_1,x_3]=x_5, [x_1,x_5]=x_6, 
[x_2,x_4]=x_6, [x_3,x_5]=fx_6.$$
However, for $f\neq 0$, this algebra is isomorphic to $L_{6,19}(f)$
by $x_1\mapsto y_1+(1/f)y_3$, and $x_i\mapsto y_i$ for $i>1$ (where 
the $x_i$ are the basis elements of $K^1_f$). A Gr\"{o}bner basis
computation reveals that $K^1_0\not\cong L_{6,19}(\epsilon)$ 
(for no $\epsilon$). Hence we get a new Lie algebra
$$L_{6,20}:  [x_1,x_2]=x_4, [x_1,x_3]=x_5, [x_1,x_5]=x_6, 
[x_2,x_4]=x_6. $$ 
\end{enumerate}
\item $d=0$. Then $f\neq 0$ as otherwise the matrix (\ref{eqnL58}) has 
rank $<2$.
We divide to get $f=1$. And to preserve this we set $a_{11}=a_{33}=1$.
Choose $a_{52}$ so that $c\mapsto 0$.
Choose $a_{11}=1$ and $a_{22}=1/a$ so that $a\mapsto 1$. We get 
$$K^2_b:  [x_1,x_2]=x_4, [x_1,x_3]=x_5, [x_1,x_4]=x_6, [x_1,x_5]=bx_6, 
[x_3,x_5]=x_6.$$
If we multiply $x_1,\ldots ,x_6$ by $1,b^2,b,b^2,b,b^2$ respectively,
then $b$ disappears if it is nonzero. So we have two cases: $b=0,1$. 
However $K^2_1\cong K^2_0$ by $x_1\mapsto y_1+y_3$ and $x_i\mapsto y_i$ for
$i>1$. And $L_{6,20}\cong K^2_0$ by $x_1\mapsto y_1$, $x_2\mapsto y_3$,
$x_3\mapsto y_2$, $x_4\mapsto y_5$, $x_5\mapsto y_4$, $x_6\mapsto y_6$.
\end{enumerate}

Now suppose that the characteristic is $2$ and $d=f=0$. Then $e\neq 0$,
otherwise the matrix (\ref{eqnL58}) has rank $<2$. We divide to get $e=1$,
Choose $a_{32}=a_{23}=0$ and $a_{53}$ so that $c\mapsto 0$. Furthermore,
choose $a_{31}$ so that $a\mapsto 0$ and $a_{21}$ so that $b\mapsto 0$.
We get the Lie algebra
$$[x_1,x_2]=x_4, [x_1,x_3]=x_5, [x_2,x_5]=x_6, [x_3,x_4]=x_6$$
(which for characteristic not $2$ is isomorphic to $L_{6,19}(-1)$).

\subsection{$L_{5,9}$}\label{sec_6_59}

The automorphism group consists of
$$ \begin{pmatrix} a_{11} & a_{12} & 0 & 0 & 0 \\
                   a_{21} & a_{22} & 0 & 0 & 0\\
                   a_{31} & a_{32} & \delta & 0 & 0 \\
                   a_{41} & a_{42} & a_{11}a_{32}-a_{31}a_{12} & a_{11}\delta
 & a_{12}\delta \\
                   a_{51} & a_{52} & a_{21}a_{32}-a_{31}a_{22} & a_{21}\delta 
& a_{22}\delta
\end{pmatrix},$$
where $\delta = a_{11}a_{22}-a_{12}a_{21}$.

$H^2(L_{5,9},F)$ consists of $\theta = a\Delta_{14}+b(\Delta_{15}+\Delta_{24})
+c\Delta_{25}$. Moreover, $\theta^{\perp}\cap C(L_{5,9})=0$ if and
only if $ac-b^2\neq 0$. The automorphism group acts as follows
\begin{align*}
a &\mapsto (a_{11}^2a+2a_{11}a_{21}b+a_{21}^2c)\delta\\
b &\mapsto (a_{11}a_{12}a+(a_{11}a_{22}+a_{21}a_{12})b+a_{21}a_{22}c)\delta\\
c &\mapsto (a_{12}^2a+2a_{12}a_{22}b+a_{22}^2c)\delta.
\end{align*}

We claim that, unless $a=c=0$ and the characteristic is $2$, the $a_{ij}$
can be chosen so that $b$ is mapped to $0$. Note that $b\mapsto
\delta( a_{12}(a_{11}a+a_{21}b)+a_{22}(a_{11}b+a_{21}c))$. If
$c\neq 0$ then choose $a_{12}=0$ and $a_{21}=-a_{11}b/c$. If $c=0$, but
$a\neq 0$, then choose $a_{21}=0$ and $a_{12}=-a_{22}b/a$. If $a=c=0$, then
choose the $a_{ij}$ so that $a_{11}a_{22}-a_{21}a_{12}\neq 0$ and
$a_{11}a_{22}+a_{12}a_{21}=0$. This can be
done if the characteristic is not $2$. \par
So if $\chi\neq 2$ or one of $a,c$ is nonzero, then
we may assume that $b=0$. Since $ac-b^2\neq 0$ we have that
$a,c\neq 0$. Divide to get $a=1$, yielding
$$L_{6,21}(\epsilon): [x_1,x_2]=x_3, [x_1,x_3]=x_4, [x_2,x_3]=x_5, 
[x_1,x_4]=x_6, [x_2,x_5]= \epsilon x_6.$$
(Here we write $\epsilon$ in place of $c$. We have $\epsilon\neq 0$, to have 
a $1$-dimensional centre.)
By multiplying $x_1,\ldots,x_6$ by $1,\alpha,\alpha,\alpha,\alpha^2,
\alpha$ we see that $L_{6,21}(\epsilon)\cong L_{6,21}(\alpha^2\epsilon)$.
If there is an isomorphism $\sigma : L_{6,21}(\epsilon) \to L_{6,21}(\delta)$
then by Lemma \ref{lem1.4} there is an automorphism $\phi$ of $L_{5,9}$
such that $\phi(\Delta_{14}+\epsilon\Delta_{25}) =\alpha (\Delta_{14}
+\delta \Delta_{15})$. This leads to the polynomial equations
$a_{11}a_{12}+a_{21}a_{22}\epsilon=0$, $a_{12}^2+a_{22}^2\epsilon = (a_{11}^2 +
a_{21}^2\epsilon)\delta$. We add the equation $D_1(a_{11}a_{22}-a_{12}
a_{21})-1=0$, $D_2(a_{11}^2+a_{21}^2\epsilon)-1=0$
and compute a Gr\"{o}bner basis. It contains the polynomials
$a_{11}^2\delta - a_{22}^2\epsilon$, $a_{12}^2 - a_{21}^2\epsilon\delta$. 
By Lemma \ref{lem3.2},
it follows that $\epsilon=\beta^2\delta$ for some $\beta\in F^*$.
We conclude that $L_{6,21}(\epsilon)\cong L_{6,21}(\delta)$ if and only if 
there is a $\beta\in F^*$ such that $\epsilon=\beta^2\delta$.\par
If the characteristic is $2$, and $a=c=0$, then we divide by $b$ and get
the Lie algebra
$$[x_1,x_2]=x_3, [x_1,x_3]=x_4, [x_2,x_3]=x_5, [x_1,x_5]=x_6,
[x_2,x_4]=x_6.$$

\subsection{$L_{4,1}$}

Here we find the $2$-dimensional central extensions of the abelian Lie
algebra $L_{4,1}$ of dimension $4$. In this case $H^2(L_{4,1},F)$ consists of 
all skew-symetric $\theta : L_{4,1}\times L_{4,1}\to F$. 
The automorphism group of 
$L_{4,1}$ is $\GL_4(F)$. We need to classify the orbits of this group
on $2$-dimensional subspaces of $H^2(L_{4,1},F)$. We denote the basis elements
of such a subspace by $\theta_1$ and $\theta_2$. By Theorem \ref{diagthm}
we may assume that $\theta_1=\Delta_{12}+\Delta_{34}$, or
$\theta_1=\Delta_{12}$. \par
We start with the first case, $\theta_1=\Delta_{12}+\Delta_{34}$. Set
$$H=\left\{ \begin{pmatrix} a_{11} & a_{12} & 0 & 0\\
             a_{21} & a_{22} & 0 & 0\\
              0 & 0 & a_{33} & a_{34}\\
              0 & 0 & a_{43} & a_{44} \end{pmatrix}
\mid a_{11}a_{22}-a_{12}a_{21} = 
a_{33}a_{44}-a_{34}a_{43} = u \neq 0 \right\}.$$
Then $H$ is a subgroup of $\Aut(L_{4,1})$ such that $\phi\cdot \theta_1 =
u\theta_1$. Hence $\phi\cdot \langle \theta_1,\theta_2 \rangle =
\langle \theta_1, \phi\cdot \theta_2\rangle$, so that we only have
to consider the action of $H$ on $\theta_2$. After, if necessary,
subtracting a multiple of $\theta_1$, we may assume that
$\theta_2 = a\Delta_{13}+b\Delta_{14}+c\Delta_{23}
+d\Delta_{24}+e\Delta_{34}$.
Note that $\phi\cdot \theta_2(x_1,x_2) =u\theta_2(x_1,x_2)=0$, and
hence also $\phi\cdot \theta_2$ has this form. Let $\phi\in H$. Then
$\phi$ acts as follows on $\theta_2$:

\begin{align*}
a & \mapsto a_{11}a_{33}a + a_{11}a_{43}b
+a_{21}a_{33}c +a_{21}a_{43}d \\
b&\mapsto a_{11}a_{34}a +a_{11}a_{44}b
+ a_{21}a_{43}c +a_{21}a_{44} d\\
c&\mapsto a_{12}a_{33}a +a_{12}a_{43}b
+a_{22}a_{33}c+a_{22}a_{43}d\\
d&\mapsto a_{12}a_{34}a+a_{12}a_{44}b
+a_{22}a_{34}c+a_{22}a_{44}d\\
e &\mapsto ue.
\end{align*}

Suppose that $a\neq 0$, then after dividing we have $a=1$.
In order to preserve this we set $a_{21}=a_{43}=0$,
and $a_{11}a_{33}=1$. Then we choose $a_{34}$ and
$a_{12}$ such that $b$, $c$ are mapped to $0$.
By choosing $u$ we can get $e=0,1$. This leads to two
Lie algebras
$$ L_{6,22}(\epsilon) : [x_1,x_2]=x_5, [x_1,x_3]=x_6, [x_2,x_4]=\epsilon x_6,
[x_3,x_4]=x_5$$
(corresponding to $e=0$, where we write $\epsilon$ in place of $d$), 
$$ K_{\epsilon}^1 : [x_1,x_2]=x_5, [x_1,x_3]=x_6, [x_2,x_4]=\epsilon x_6,
[x_3,x_4]=x_5+x_6$$
(corresponding to $e=1$). If the characteristic is not $2$, then
for all $\epsilon$ there is a $\delta$ such that $K_{\epsilon}^1\cong 
L_{6,22}(\delta)$. If $\epsilon\neq 0$ and $\epsilon\neq -\frac{1}{4}$ then 
we set $\delta=\epsilon+\frac{1}{4}$, and define $\psi$
by $\psi(x_1) = (\frac{1}{4\epsilon}+1)y_1+\frac{1}{2\epsilon}y_4$,
$\psi(x_2)=\epsilon y_3$, $\psi(x_3)=y_2+\frac{1}{2}y_3$, $\psi(x_4)=y_4$,
$\psi(x_5)=-\frac{1}{2}y_5 +(\epsilon+\frac{1}{4})y_6$, $\psi(x_6)=y_5$.
(Here the $x_i$ are the basis elements of $K_{\epsilon}^1$.)
If $\epsilon=0$, then we set $\delta=1$ and $\psi(x_1)=2y_4$, $\psi(x_2)=
-\frac{1}{4} y_2+\frac{1}{4}y_3$, $\psi(x_3)=-\frac{1}{2}y_3$, $\psi(x_4) = 
y_1-y_4$, $\psi(x_5)=-\frac{1}{2}y_5 + \frac{1}{2} y_6$, $\psi(x_6) = 
y_5$. If $\epsilon=-\frac{1}{4}$, then we set $\delta=0$ and $\psi(x_1) = 
-2y_1-4y_4$, $\psi(x_2)=\frac{1}{2}y_2 +\frac{1}{2}y_3$, 
$\psi(x_3)=-2y_2-y_3$, $\psi(x_4)=y_1+y_4$, $\psi(x_5)=y_5-y_6$,
$\psi(x_6)=2y_6$. We conclude that if the characteristic is not $2$, then
we can discard $K_{\epsilon}^1$. (However, if the characteristic is 2, then 
a Gr\"{o}bner basis computation shows that they are not isomorphic.)\par
Now we consider the question whether $L_{6,22}(\epsilon)\cong L_{6,22}
(\delta)$. With $\gamma =1/\beta \in F^*$ and $y_1=\gamma x_1$, $y_2=x_2$, 
$y_3=\gamma y_3$, $y_4=x_4$, $y_5=\gamma y_5$, $y_6=\gamma^2 x_6$ we see that
$L_{6,22}(\epsilon) \cong L_{6,22}(\beta^2\epsilon)$. However by a 
Gr\"{o}bner basis computation (which works if $\chi\neq 2$), we
have that isomorphism implies that $x^2\delta-y^2\epsilon=0$, $u^2-v^2\epsilon
\delta=0$ has a solution with $xy-uv\neq 0$. Hence by Lemma \ref{lem3.2} 
there is a $\beta\in F^*$ with $\delta=\beta^2 \epsilon$. 
(If $\chi=2$ and the field is perfect, then $\delta/\epsilon$ is always a 
square, and hence $L_{6,22}(\epsilon)\cong L_{6,22}(\delta)$.) 
Conclusion: $L_{6,22}(\epsilon)\cong L_{6,22}(\delta)$ if and only if there 
is a $\beta\in F^*$ such that $\delta=\beta^2 \epsilon$. \par
Now suppose that $a=0$. Then if not all three of $b,c,d$ are $0$, we
can make $a$ nonzero, and are back in the previous case. 
If $a=b=c=d=0$, then $\theta_2 = \Delta_{34}$.
After subtracting it from $\theta_1$ we get that $\theta_1=\Delta_{12}$.
This yields
$$K^2: [x_1,x_2] = x_5, [x_3,x_4]=x_6.$$
However, $K^2\cong K_0^1$ by $x_1\mapsto y_4$, $x_2\mapsto -y_2-y_3$, 
$x_3\mapsto y_2$, $x_4\mapsto -y_1-y_4$, $x_5\mapsto y_5+y_6$, 
$x_6\mapsto y_5$. \par
Now suppose that $\theta_1=\Delta_{12}$. As before we may assume that
$$\theta_2 = a\Delta_{13}+b\Delta_{14}+c\Delta_{23}
+d\Delta_{24}+e\Delta_{34}.$$
The space $\theta_1^\perp$ is spanned by $x_3,x_4$. A small calculation 
shows that $\theta_1^\perp \cap\theta_2^\perp = 0$ if and only if
the matrix
\begin{equation}\label{eqn_4_ab}
\begin{pmatrix}
a & b \\
c & d \\
0 & e \\
-e & 0 
\end{pmatrix}
\end{equation}
has rank $2$. Set $\theta = x\theta_1 +\theta_2$. Then $\theta^\perp =0$
if and only if $xe-ad+bc \neq 0$. If $ad-bc\neq 0$ then we choose $x=0$.
If $ad-bc =0$, then since the matrix (\ref{eqn_4_ab}) has rank $2$,
$e\neq 0$. In this case we choose $x=1$. The conclusion is that $x$
can be chosen such that $\theta^{\perp}=0$. By Theorem \ref{diagthm}
there is a basis of $L_{4,1}$ such that $\theta= \Delta_{12}+\Delta_{34}$.
We set $\theta_1=\theta$, and we are back in the previous case.

\subsection{$L_{4,2}$}\label{sec6_42}

In this section we classify $2$-dimensional central extensions
of the $4$-dimensional Lie algebra $L_{4,2}$. In Section \ref{secnonab4_1}
it was shown that there are three orbits of $\Aut(L_{4,2})$ on 
$H^2(L_{4,2},F)$. Here we have to classify $2$-dimensional subspaces of 
$H^2(L_{4,2},F)$, and we may assume that the first basis element, $\theta_1$
is one of the three computed in Section \ref{secnonab4_1}. \par
Suppose that $\theta_1= \Delta_{13}+\Delta_{24}$. Let $\phi$ be an
automorphism of $L_{4,2}$ as considered in section \ref{secnonab4_1}. A small
calculation shows that $\phi\theta_1$ is a scalar multiple of
$\theta_1$ if and only if $a_{12}=0$, $a_{34}=-a_{11}a_{21}$, 
$a_{44}=a_{11}^2$. After maybe adding a scalar multiple of $\theta_1$, 
we may assume that $\theta_2=b\Delta_{14}+c\Delta_{23}+d\Delta_{24}$. 
With $\phi$ satisfying $a_{12}=0$, $a_{34}=-a_{11}a_{21}$, $a_{44}=a_{11}^2$ 
we have
\begin{align*}
\phi\theta_2(x_1,x_3) &= a_{11}a_{22}a_{21}c\\
\phi\theta_2(x_1,x_4) &= a_{11}(a_{11}^2b-a_{21}^2c+a_{11}a_{21}d)\\
\phi\theta_2(x_2,x_3) &= a_{11}a_{22}^2 c\\
\phi\theta_2(x_2,x_4) &= a_{11}a_{22}(-a_{21}c+a_{11}d)
\end{align*}

We consider some case distinctions:
\begin{enumerate}
\item $c=0$ and $d\neq 0$. Then we divide by $d$, and use $a_{21}$
to get $b=0$. Then $\theta_2=\Delta_{24}$. After subtracting it
from $\theta_1$ we get $\theta_1=\Delta_{13}$. This gives the Lie algebra
$$K^1: [x_1,x_2]=x_3, [x_1,x_3]=x_5, [x_2,x_4]= x_6.$$
However, this Lie algebra is isomorphic to $L_{6,19}(0)$ by
$x_1\mapsto y_2$, $x_2\mapsto y_1$, $x_3\mapsto -y_4$, $x_4\mapsto
y_3$, $x_5\mapsto -y_6$, $x_6\mapsto y_5$. 
\item $c=d=0$. Then we divide and get $b=1$ and $\theta_2=\Delta_{14}$.
This leads to the Lie algebra
$$L_{6,23}: [x_1,x_2]=x_3, [x_1,x_3]=x_5, [x_1,x_4]=x_6, [x_2,x_4]= x_5.$$
A Gr\"{o}bner basis computation shows that this Lie algebra
is not isomorphic to $K^1$ (over any field). 
\item $c\neq 0$. If $d\neq 0$ and $\chi\neq 2$, then we set 
$a_{21}=1$, $a_{11}=2c/d$, $a_{22}=d/2c$. This means that $\phi\theta_2 = 
c\Delta_{13}+b'\Delta_{14}+c'\Delta_{23}+c\Delta_{24}$, with $c'\neq 0$.
 Now we subtract
$c\theta_1$, and divide by $c'$. We get $\theta_2=b\Delta_{14}+\Delta_{23}$. 
If $d=0$, then we set $a_{21}=0$ and divide by $c$ and arrive at the
same conclusion. This yields
$$ L_{6,24}(\epsilon): [x_1,x_2]=x_3, [x_1,x_3]=x_5, [x_1,x_4]=\epsilon x_6, 
[x_2,x_3]=x_6, [x_2,x_4]= x_5.$$
Because of different lower central series dimensions, this Lie algebra
is not isomorphic to $K^1$ or $L_{6,23}$ above.\par 
By setting $y_1=\alpha x_1$, $y_2=x_2$, $y_3=\alpha x_3$, 
$y_4=\alpha^2 x_4$, $y_5=\alpha^2 x_5$, $y_6=\alpha x_6$ we get that
$L_{6,24}(\epsilon)\cong L_{,24}(\alpha^2 \epsilon)$ for any $\alpha\in F^*$. 
By a Gr\"{o}bner basis computation (which works for all
fields of characteristic not $2$), we get that $L_{6,24}(\epsilon)\cong 
L_{6,24}(\delta)$
implies that $x^2-y^2\epsilon\delta=0$ and $s^2\delta-t^2\epsilon=0$ 
have a solution in $F$ with $st-xy\neq 0$. By Lemma \ref{lem3.2},
there is an $\alpha\in F^*$ such that $\delta=\alpha^2 \epsilon$. 
\end{enumerate}
Suppose that $\theta_1=\Delta_{13}$. Let $\phi\in \Aut(L_{4,2})$ be such 
that $a_{34}=a_{12}=0$. Then $\phi\theta_1 = a_{11}^2a_{22}\theta_1$.
We may assume that $\theta_2 = b\Delta_{14}+c\Delta_{23}+d\Delta_{24}$.
A small calculation shows that $\theta_1^\perp \cap \theta_2^\perp 
\cap C(L) =0$
if and only if $d\neq 0$ or $b\neq 0$. The automorphism $\phi$
with $a_{34}=a_{12}=0$ acts as follows
\begin{align*}
b&\mapsto a_{11}a_{44}b+a_{21}a_{44}d\\
c&\mapsto a_{11}a_{22}^2c\\
d&\mapsto a_{22}a_{44}d. 
\end{align*}
(The value of $\phi\theta_2(x_1,x_3)$ does not matter, because we can
make it zero afterwards by subtracting a suitable multiple of $\theta_1$.)\par
If $d\neq 0$ then after dividing it is $1$. By choosing $a_{21}$
appropriately we can get $b\mapsto 0$. With $a_{11}$ we can enforce
$c\mapsto 1$ if $c\neq 0$. So we get $\theta_2 = *\Delta_{23}+\Delta_{24}$,
where $*=0,1$. If $*=1$ then we get the Lie algebra
$$K^2 : [x_1,x_2]=x_3, [x_1,x_3]=x_5, [x_2,x_3]=x_6, [x_2,x_4]=x_6.$$
However, if $\chi\neq 2$ then $K\cong L_{6,24}(1)$ by $x_1\mapsto y_1-y_2$,
$x_2\mapsto y_1+y_2$, $x_3\mapsto 2y_3$, $x_4\mapsto y_3+y_4$,
$x_5\mapsto 2y_5-2y_6$, $x_6\mapsto 2y_5+2y_6$. If $*=0$ then
we have $\theta_1=\Delta_{13}$, $\theta_2=\Delta_{24}$, which are
also the cocycles defining $K^1$.\par
If $d=0$ then $b\neq 0$; divide to get $b=1$. By choosing $a_{11}$
we get $c\mapsto 0,1$. If $c=1$ then we get  
$$K^3 : [x_1,x_2]=x_3, [x_1,x_3]=x_5, [x_1,x_4]=x_6, [x_2,x_3]=x_6.$$
However, $K^3\cong L_{6,24}(0)$ by $x_1\mapsto -y_2$, $x_2\mapsto -y_1$,
$x_3\mapsto -y_3$, $x_4\mapsto -y_4$, $x_5\mapsto y_6$, $x_6\mapsto y_5$.
The choice $c=0$ leads to
$$L_{6,25} : [x_1,x_2]=x_3, [x_1,x_3]=x_5, [x_1,x_4]=x_6.$$
By lower central series dimensions, this Lie algebra could by isomorphic
to $K^1$ or $L_{6,23}$. However, Gr\"{o}bner basis computations show that
it is not (over any field).\par
Now suppose that $\theta_1=\Delta_{14}$. Let $\phi\in \Aut(L_{4,2})$ 
be such that $a_{12}=0$. Then $\phi\theta_1=a_{11}a_{44}\theta_1$.
We may assume that $\theta_2 = a\Delta_{13}+c\Delta_{23}+d\Delta_{24}$.
We have that $\theta_1^\perp\cap\theta_2^\perp\cap C(L) =0$ if and
only if $a\neq 0$ or $c\neq 0$. The automorphism $\phi$ acts as follows
\begin{align*}
a&\mapsto a_{11}a_{22}(a_{11}a+a_{21}c)\\
c&\mapsto a_{11}a_{22}^2c\\
d&\mapsto a_{22}a_{34}c+a_{22}a_{44}d. 
\end{align*}
If $c\neq 0$ then we divide to make it $1$. Then by choosing
$a_{34}$ and $a_{21}$ we can ensure that $d\mapsto 0$, $a\mapsto 0$
respectively. So we get $\theta_2=\Delta_{23}$ and
$$K^4 : [x_1,x_2]=x_3, [x_1,x_4]=x_5, [x_2,x_3]=x_6.$$
However, $K^4\cong K^1$ by $x_1\mapsto y_2$, $x_2\mapsto -y_1$, 
$x_3\mapsto y_3$, $x_4\mapsto -y_4$, $x_5\mapsto -y_6$, $x_6\mapsto
-y_5$. \par
If $c=0$ then $a\neq 0$. After dividing we get $a=1$. By choosing
$a_{44}$ we get $d=0,1$. If $d=1$ then we get the same cocycles that 
define $L_{6,23}$.  If $d=0$ then we get the same cocycles as were used to 
define $L_{6,25}$.

\subsection{$L_{4,3}$}

In this case $H^2(L_{4,3},F)$ is $2$-dimensional, and spanned by 
$\Delta_{14}$, $\Delta_{23}$. This leads to the Lie algebra
$$[x_1,x_2]=x_3, [x_1,x_3]=x_4, [x_1,x_4]=x_5, [x_2,x_3]=x_6.$$
By interchanging $x_5$ and $x_6$ we see that this algebra is isomorphic
to $L_{6,21}(0)$. 

\subsection{$L_{3,1}$}

Here $H^2(L_{3,1},F)$ is $3$-dimensional. We get the Lie algebra
$$L_{6,26}: [x_1,x_2]=x_4, [x_1,x_3]=x_5, [x_2,x_3]=x_6.$$

\subsection{The list}

First we get nine algebras denoted $L_{6,k}$ for $k=1,\ldots,9$ that are the 
direct sum of $L_{5,k}$ and a $1$-dimensional abelian ideal. Subsequently
we get the following Lie algebras.

\begin{enumerate}
\item[$L_{6,10}$] $[x_1,x_2]=x_3, [x_1,x_3]=x_6, [x_4,x_5]=x_6.$
\item[$L_{6,11}$] $[x_1,x_2]=x_3, [x_1,x_3]=x_4, [x_1,x_4]=x_6, [x_2,x_3]=x_6,
[x_2,x_5]=x_6.$
\item[$L_{6,12}$] $[x_1,x_2]=x_3, [x_1,x_3]=x_4, [x_1,x_4]=x_6, [x_2,x_5]=x_6.$
\item[$L_{6,13}$] $[x_1,x_2]=x_3, [x_1,x_3]=x_5, [x_2,x_4]=x_5, [x_1,x_5]=x_6,
[x_3,x_4]=x_6.$
\item[$L_{6,14}$] $[x_1,x_2]=x_3, [x_1,x_3]=x_4, [x_1,x_4]=x_5, [x_2,x_3]=x_5, 
[x_2,x_5]=x_6,[x_3,x_4]=-x_6.$
\item[$L_{6,15}$] $[x_1,x_2]=x_3, [x_1,x_3]=x_4, [x_1,x_4]=x_5, [x_2,x_3]=x_5, 
[x_1,x_5]=x_6,[x_2,x_4]=x_6.$
\item[$L_{6,16}$] $[x_1,x_2]=x_3, [x_1,x_3]=x_4, [x_1,x_4]=x_5, [x_2,x_5]=x_6,
[x_3,x_4]=-x_6.$
\item[$L_{6,17}$]  $[x_1,x_2]=x_3, [x_1,x_3]=x_4, [x_1,x_4]=x_5, [x_1,x_5]=x_6,
[x_2,x_3]= x_6.$
\item[$L_{6,18}$]  $[x_1,x_2]=x_3, [x_1,x_3]=x_4, [x_1,x_4]=x_5, 
[x_1,x_5]=x_6.$
\item[$L_{6,19}(\epsilon)$] $[x_1,x_2]=x_4, [x_1,x_3]=x_5, [x_2,x_4]=x_6, 
[x_3,x_5]=\epsilon x_6.$
Isomorphism: $L_{6,19}(\epsilon)\cong L_{6,19}(\delta)$ if and only if 
there is an 
$\alpha\in F^*$ such that $\delta=\alpha^2\epsilon$. 
\item[$L_{6,20}$] $[x_1,x_2]=x_4, [x_1,x_3]=x_5, [x_1,x_5]=x_6, 
[x_2,x_4]=x_6. $
\item[$L_{6,21}(\epsilon)$] $[x_1,x_2]=x_3, [x_1,x_3]=x_4, [x_2,x_3]=x_5, 
[x_1,x_4]=x_6, [x_2,x_5]= \epsilon x_6.$
Isomorphism: $L_{6,21}(\epsilon)\cong L_{6,21}(\delta)$ if and only if there 
is an $\alpha\in F^*$ such that $\delta=\alpha^2\epsilon$. 
\item[$L_{6,22}(\epsilon)$] $[x_1,x_2]=x_5, [x_1,x_3]=x_6, [x_2,x_4]=
\epsilon x_6, [x_3,x_4]=x_5.$ 
Isomorphism: $L_{6,22}(\epsilon)\cong L_{6,22}(\delta)$ if and only if there 
is an $\alpha\in F^*$ such that $\delta=\alpha^2\epsilon$.
\item[$L_{6,23}$] $[x_1,x_2]=x_3, [x_1,x_3]=x_5, [x_1,x_4]=x_6, 
[x_2,x_4]= x_5.$
\item[$L_{6,24}(\epsilon)$] $[x_1,x_2]=x_3, [x_1,x_3]=x_5, [x_1,x_4]=\epsilon 
x_6, [x_2,x_3]=x_6, [x_2,x_4]= x_5.$
Isomorphism: $L_{6,24}(\epsilon)\cong L_{6,24}(\delta)$ if and only if there 
is an $\alpha\in F^*$ such that $\delta=\alpha^2\epsilon$.
\item[$L_{6,25}$] $[x_1,x_2]=x_3, [x_1,x_3]=x_5, [x_1,x_4]=x_6.$
\item[$L_{6,26}$] $[x_1,x_2]=x_4, [x_1,x_3]=x_5, [x_2,x_3]=x_6.$
\end{enumerate}

\section{Example}\label{rec_exa}

In this section we show in an example how the recognition procedure works.
Let $L$ be the Lie algebra with basis ${x_1,\ldots,x_6}$ and 
$$[x_1,x_2]=x_3+x_6, [x_1,x_3]=x_5+x_6, [x_1,x_4]=x_5, [x_2,x_3]=x_6, 
[x_2,x_4]=2x_5+x_6.$$

Then $C(L)$ is spanned by $x_5,x_6$ and $L/C(L)$ is isomorphic to 
$L_{4,2}$. The cocycles that define $L$ as a central extension of
$L_{4,2}$ are $\theta_1= \Delta_{13}+\Delta_{14}+2\Delta_{24}$,
and $\theta_2=\Delta_{12}+\Delta_{13}+\Delta_{23}+\Delta_{24}$. 
We denote the basis elements of $L_{4,2}$ by $y_1,\ldots,y_4$.
We combine the operations of changing the basis of $L_{4,2}$, and
taking linear combinations of $\theta_1$, $\theta_2$. \par
First we change the basis of $L_{4,2}$ such that $\theta_1$
has one of the three forms described in Section \ref{secnonab4_1}.
For that note that $a=b=1$ and $d=2$ (notation as in Section
\ref{secnonab4_1}). Hence we choose $a_{44}=\frac{1}{2}$ and 
$a_{34}=-\frac{1}{2}$. This means that the new basis of 
$L_{4,2}$ is $u_1,\ldots,u_4$, where $u_4=-\frac{1}{2}y_3+\frac{1}{2}y_4$
and $u_i=y_i$ otherwise. Relative to this basis $\theta_1=
\Delta_{13}+\Delta_{24}$ and $\theta_2 = \Delta_{12}+\Delta_{13}-
\frac{1}{2}\Delta_{14}+\Delta_{23}$. Now we transform $\theta_2$ as
described in Section \ref{sec6_42}. First we subtract $\theta_1$
to get rid of the $\Delta_{13}$. So we set $\theta_2'=
-\theta_1+\theta_2 = \Delta_{12}-\frac{1}{2}\Delta_{14}+\Delta_{23}-
\Delta_{24}$. In the notation of Section \ref{sec6_42} we have
$b=-\frac{1}{2}$, $c=1$, $d=-1$. So we set $a_{11}=-2$, 
$a_{21}=1$, $a_{22}=-\frac{1}{2}$. Because of the condition on
the automorphism (that is imposed to map $\theta_1$ to a scalar
multiple of itself), we have $a_{34}=2$, $a_{44}=4$. So we get
a new basis of $L_{4,2}$, given by $v_1=-2u_1+u_2$, $v_2=
-\frac{1}{2} u_2$, $v_3=u_3$, $v_4=2u_3+4u_4$. Relative to this
basis $\theta_1 = -2\Delta_{13}-2\Delta_{24}$, so we set 
$\theta_1'=-\frac{1}{2} \theta_1$. Relative to the new basis
$\theta_2' = \Delta_{12}+\Delta_{13}+2\Delta_{14}-\frac{1}{2}
\Delta_{23}+\Delta_{24}$. To get rid of the $\Delta_{13}$ and
$\Delta_{24}$ we subtract $\theta_1'$, i.e., we set 
$\theta_2'' = -\theta_1'+\theta_2'' = -\frac{1}{2} \theta_1+\theta_2$.
Then $\theta'' = \Delta_{12}+2\Delta_{14}-\frac{1}{2}\Delta_{23}$.
According to the procedure outlined in Section \ref{sec6_42} we
now have to divide by the coefficient of $\Delta_{23}$, i.e, we
set $\theta_2'''=-2\theta_2'' = \theta_1-2\theta_2$. 
Hence $\theta_2'''= -2\Delta_{12}-4\Delta_{14}+\Delta_{23}$.\par
Now let $K$ be the Lie algebra that is the central extension of
$L_{4,2}$ with respect to $\theta_1'$ and $\theta_2'''$. Denoting
the basis elements of $K$ by $w_i$ we get the multiplication table
$$[w_1,w_2]=w_3+w_4, [w_1,w_3]=w_5, [w_1,w_4]=-4w_6, [w_2,w_3]=w_6,
[w_2,w_4]=w_5.$$
In order to construct an isomorphism between $L$ and $K$ we map
$v_i$ to $w_i$ for $1\leq i\leq 4$. This means that $x_1\mapsto 
-\frac{1}{2}w_1-w_2$, $x_2\mapsto -2w_2$, $x_3\mapsto w_3$,
$x_4\mapsto \frac{1}{2}v_4$. The changes from $\theta_1$,
$\theta_2$ to $\theta_1'$, $\theta_2'''$ correspond to a base
change in the centre. According to the proof of Lemma \ref{lem1.4},
the matrix of the base change in the centre is the transpose of 
the matrix describing the change in the $\theta$'s. This means that
$x_5\mapsto -\frac{1}{2} w_5+w_6$ and $x_6\mapsto -2w_6$. \par
Finally, in order to construct an isomorphism with a Lie algebra
from the list we have to subtract a coboundary (the $\Delta_{12}$
occurring in $\theta_2'''$). Then we get an isomorphism $L\to L_{6,24}(-4)$
by $w_3\mapsto z_3+2z_6$, and $w_i\mapsto z_i$ for $i\neq 3$
(where the $z_i$ are the basis elements of $L_{6,24}(-4)$).

\section{Comments}\label{sec_comm}

With the help of the implementation of recognition procedure it is
now straightforward to compare other classifications of $6$-dimensional
nilpotent Lie algebras. Because the programs will be available in
{\sf GAP} and {\sf Magma}, we do not include a table linking the
Lie algebras given here to the ones in various other papers. Instead
we briefly give the conclusions that we obtained using the programs.\par 
In \cite{morozov} V. Morozov has obtained a classification of
$6$-dimensional nilpotent Lie algebras over a field of characteristic $0$.
His classification is essentially the same as ours. However, for two of the 
parametric classes it seems to lack the conditions for isomorphism given 
here. The classes with numbers 14 and 18 in \cite{morozov} depend on a
parameter $\gamma$, without restrictions on $\gamma$. They are isomorphic to 
$L_{6,19}(\gamma)$ and $L_{6,21}(\gamma)$ respectively. \par
The classifications by R. Beck and B. Kolman (\cite{beko3}) and O. Nielsen
(\cite{nielsen}) of $6$-dimensional nilpotent Lie algebras over $\R$ are the 
same as the one given here. The classification by M.-P. Gong (over 
algebraically closed fields) lacks a Lie algebra isomorphic to 
$L_{6,22}(\epsilon)$ with $\epsilon\neq 0$.
\par
C. Schneider (\cite{schneider2}) obtains 34 nilpotent Lie algebras of 
dimension $6$, over the fields $\F_3$ and $\F_5$. This corresponds to 
our results. Furthermore, in the context of the {\sf liealgdb} package
(\cite{liealgdb}) for {\sf GAP} he has obtained a list of $34$ nilpotent 
Lie algebras that is valid over any finite field of characteristic not $2$. 
With the recogonition procedure we have checked that his classification 
agrees with the one given here.

\def\cprime{$'$} \def\cprime{$'$} \def\cprime{$'$}


\begin{thebibliography}{10}

\bibitem{beko3}
R.~E. Beck and B.~Kolman.
\newblock Construction of nilpotent {L}ie algebras over arbitrary fields.
\newblock In Paul~S. Wang, editor, {\em Proceedings of the 1981 ACM Symposium
  on Symbolic and Algebraic Computation}, pages 169--174. ACM New York, 1981.

\bibitem{Magma}
W.~Bosma, J.~Cannon, and C.~Playoust.
\newblock The {M}agma algebra system. {I}. {T}he user language.
\newblock {\em J. Symbolic Comput.}, 24(3-4):235--265, 1997.
\newblock Computational algebra and number theory (London, 1993).

\bibitem{liealgdb}
M.~Costantini, W.~A. de~Graaf, and C.~Schneider.
\newblock {\sf liealgdb}, a database of {L}ie algebras.
\newblock a {\sf GAP}4 package.
\newblock in preparation.

\bibitem{gra11}
W.~A. de~Graaf.
\newblock Classification of solvable {L}ie algebras.
\newblock {\em Experiment. Math.}, 14(1):15--25, 2005.

\bibitem{GAP4}
The GAP~Group.
\newblock {\em {GAP -- Groups, Algorithms, and Programming, Version 4.4}},
  2004.
\newblock \verb+(http://www.gap-system.org)+.

\bibitem{gong}
M.-P. Gong.
\newblock {\em Classification of Nilpotent Lie Algebras of Dimension 7}.
\newblock PhD thesis, University of Waterloo, Waterloo, Canada, 1998.

\bibitem{jac2}
Nathan Jacobson.
\newblock {\em Lectures in abstract algebra}.
\newblock Springer-Verlag, New York, 1975.
\newblock Volume II: Linear algebra, Reprint of the 1953 edition [Van Nostrand,
  Toronto, Ont.], Graduate Texts in Mathematics, No. 31.

\bibitem{morozov}
V.~V. Morozov.
\newblock Classification of nilpotent {L}ie algebras of sixth order.
\newblock {\em Izv. Vys\v s. U\v cebn. Zaved. Matematika}, 1958(4
  (5)):161--171, 1958.

\bibitem{nielsen}
O.~A. Nielsen.
\newblock {\em Unitary representations and coadjoint orbits of low-dimensional
  nilpotent {L}ie groups}, volume~63 of {\em Queen's Papers in Pure and Applied
  Mathematics}.
\newblock Queen's University, Kingston, ON, 1983.

\bibitem{schneider2}
C.~Schneider.
\newblock A computer-based approach to the classification of nilpotent {L}ie
  algebras.
\newblock {\em Experiment. Math.}, 14(2):153--160, 2005.

\bibitem{skjelsund}
Tor Skjelbred and Terje Sund.
\newblock Sur la classification des alg\`ebres de {L}ie nilpotentes.
\newblock {\em C. R. Acad. Sci. Paris S\'er. A-B}, 286(5), 1978.

\end{thebibliography}
\end{document}